\documentclass[10pt]{article}
\usepackage{amsfonts,amssymb,amsmath,amscd,amsthm}
\usepackage{epsfig}

\makeatletter

\renewcommand\section{\@startsection {section}{1}{\z@}
{-30pt \@plus -1ex \@minus -.2ex}
{2.3ex \@plus.2ex}
{\normalfont\normalsize\bfseries\boldmath}}

\renewcommand\subsection{\@startsection{subsection}{2}{\z@}
{-3.25ex\@plus -1ex \@minus -.2ex}
{1.5ex \@plus .2ex}
{\normalfont\normalsize\bfseries\boldmath}}

\renewcommand{\@seccntformat}[1]{\csname the#1\endcsname. }

\makeatother

\newtheorem{theorem}{Theorem}
\newtheorem{lemma}{Lemma}

\newtheorem{corollary}{Corollary}

\theoremstyle{definition}

\newcommand{\pgap}   {{\mathcal G}}
\newcommand {\gap}     {\makebox[0.075 in]{}}   
   
\newcommand {\st}      {\gap : \gap}

\newcommand {\fto}     {\longrightarrow}

\newcommand {\set}[1]  {\left\{ {#1} \right\}}   
  
\newcommand {\pml}[1]  {{#1}^{\#}}

\begin{document}

\begin{center}
\uppercase{\bf \boldmath Eratosthenes sieve supports the $k$-tuple conjecture}
\vskip 20pt
{\bf Fred B. Holt}\\
{\tt fbholt@primegaps.info}
\end{center}
\vskip 20pt
\vskip 30pt

\centerline{\bf Abstract}
\noindent
We study Eratosthenes sieve as a discrete dynamic system.  At each stage of the sieve there is a cycle of gaps $\pgap(\pml{p})$ of
length $\phi(p^\#)$ and span $p^\#$.  There is a recursion $\pgap(p_k^\#)\fto \pgap(p_{k+1}^\#)$ that creates the next cycle from
the current one.  Under this recursion we can identify gaps and constellations among the $p$-rough numbers and the driving terms
for these gaps and constellations.  

Taking a close look at the action of the recursion on $s$ and its driving terms, we show that every admissible instance of every
admissible constellation is produced by Eratosthenes sieve.

If we take initial conditions from the cycle $\pgap(p_0^\#)$, then for all constellations of span $|s| < 2p_1$, including gaps $g < 2p_1$,
the driving terms of various lengths form Markov chains.  These yield {\it exact} models for the populations $n_s(p_k^\#)$ for all further
stages of the sieve.  These populations are all superexponential, so we factor out the superexonential growth to obtain the
exact model for the relative population $w_s(p_k^\#)$ of the constellation $s$ across all further stages of the sieve.
$$ w_{s,J}(p_k^\#) \; = \; n_{s,J}(p_k^\#) \, / \, \prod_{J+1 < p \le p_k} (p-J-1) $$

If $s$ is an admissible constellation of length $J$, then $n_{s,J}(p^\#) = \Theta \left( \prod (q-J-1)\right)$, and the asymptotic value
for the relative population is a constant $w_{s,J}(\infty) \ge 1$ that depends only on the odd prime factors that divide a span
within $s$.
\pagestyle{myheadings}
\thispagestyle{empty}
\baselineskip=12.875pt
\vskip 30pt

\section{Introduction}
We have been studying Eratosthenes sieve as a discrete dynamic system.
At each stage of the sieve, there is a cycle of gaps $\pgap(\pml{p})$
of length $\phi(\pml{p})$ (number of gaps in the cycle) and span $\pml{p}$ (sum of the gaps in the cycle).  The existence of these cycles is
well-known and the early cycles, e.g. $\pgap(\pml{5})$ or $\pgap(\pml{7})$, are routinely rediscovered.

What is novel about our approach is that we have identified a 3-step recursion that produces the next cycle of gaps from the current one,
$$\pgap(p_k^\#) \fto \pgap(p_{k+1}^\#)$$
The cycles of gaps under this recursion constitute a discrete dynamic system, and we are able to obtain analytic results far beyond what could 
be computed explicitly.

Here we show how the recursion across cycles of gaps $\pgap(\pml{p})$ 
supports the $k$-tuple conjecture.  Specifically we prove the following.

\begin{theorem}\label{ThmAdmS}
Let $s$ be a constellation of length $J$, admissible for all primes $p \le p_k$.
Then every admissible instance of $s$ modulo $p$ over all primes $p \le p_k$ has a unique occurrence in $\pgap(p_k^\#)$.
\end{theorem}

As a corollary, this provides counts of the relative populations of gaps and populations that are consistent with the estimates conjectured
by Hardy and Littlewood \cite{HL}.  The population of every admissible constellation of length $J$ grows as 
${\Theta(\prod_{p > J+1}(p-J-1))}$.

\begin{corollary}\label{CorNsj}
Let $s$ be an admissible constellation of length $J$, and let $Q$ be the product of all the odd primes that divide a span between 
boundary fusions in $s$.  Then the aggregate population of $s$ and its driving terms in $\pgap(\pml{p})$ is
$$
\sum_{j \ge J} n_{s,j} (\pml{p})  =  \prod_{q \le p} (q - \nu_q(s)) 
$$
and the asymptotic population is
\begin{eqnarray*}
n_{s,J}(\infty) & = & w_{s,J}(\infty) \cdot \prod_{q > J+1} (q-J-1)  \\
{ with} \hspace{0.2in} w_{s,J}(\infty) & = & \prod_{q \le J+1} (q-\nu_q(s))  \cdot  \prod_{\substack{q > J+1, \\ q \mid Q(s)}} \frac{q-\nu_q(s)}{q-J-1} 
\end{eqnarray*}
\end{corollary}

The context and proofs for these results are provided below.

What this tells us is that the populations of {\em all} admissible constellations of length $J$ ultimately grow at the same rate, 
$\prod (q-J-1)$.  For each admissible constellation $s$, the asymptotic relative population $w_{s,J}(\infty) \ge 1$ is a constant, which
depends only on the odd prime factors $Q(s)$ of the spans between the boundary fusions in $s$.
For any two admissible constellations $s_1, s_2$ of length $J$, the ratio of their populations is a positive constant
$$n_{s_1,J}(\pml{p}) / n_{s_2,J}(\pml{p}) \quad \xrightarrow[p \rightarrow \infty]~\quad w_{s_1,J}(\infty) / w_{s_2,J}(\infty).$$

If $s_1$ is of length $J_1$ and $s_2$ of length $J_2$, with ${J_1 < J_2}$, then
\begin{equation} 
n_{s_2,J_2}(\pml{p}) / n_{s_1,J_1}(\pml{p}) \fto 0. \label{ZeroEq}
\end{equation}
This  holds no matter how populous the longer $s_2$ is, or how rare the shorter $s_1$ is, at any stage of the sieve.

These asymptotic results evolve over incredibly large scales, far beyond the computational horizon.  Analyses of statistical samples such as
\cite{OS} must be considered in the context of this evolution \cite{FBHbias}.  Otherwise they can be misinterpreted at face value.

\section{The dynamic system}
At the stage of Eratosthenes sieve in which $p$ has been confirmed as a prime and the multiples of $p$ removed, there is a cycle of gaps
$\pgap(\pml{p})$ among the remaining candidate primes.  These remaining candidate primes are also the generators of 
$\mathbb{Z} \bmod \pml{p}$ and are also known as the $p$-rough numbers.

For example, the initial cycles of gaps are
\begin{eqnarray*}
\pgap(\pml{3}) & = & 4 \; 2 \\
\pgap(\pml{5}) & = & 6 \; 4 \; 2 \; 4 \; 2 \; 4 \; 6 \; 2
\end{eqnarray*}
The cycle $\pgap(\pml{p})$ has length $\phi(\pml{p})$ and span $\pml{p}$.

Polignac \cite{Pol} observed these cycles in 1849, as the context for his conjecture about gaps between primes.  
The smaller cycles are the basis for the wheel factoring method.  Still, it seems that every
few years the cycles $\pgap(\pml{5})$ and $\pgap(\pml{7})$ are rediscovered.

On the other hand, we have not seen any prior exposition about the recursion \cite{FBHSFU, FBHPatterns} that creates the next cycle of gaps
$\pgap(\pml{p_{k+1}})$ from the current one $\pgap(\pml{p_k})$.

\begin{lemma}
For the cycles of gaps, there is a 3-step recursion that produces $\pgap(\pml{p_{k+1}})$ from 
${\pgap(\pml{p_k}) = g_1 g_2 \ldots g_{\phi(\pml{p_k})}}$.
\begin{itemize}
\item[R1.] $p_{k+1} = g_1 + 1$
\item[R2.]  Concatenate $p_{k+1}$ copies of $\pgap(\pml{p_k})$.
\item[R3.] {\em Fusions.}  Add together $g_1+g_2$ and thereafter at the running sums indicated by the element-wise product
$p_{k+1} \ast \pgap(\pml{p_k})$.
\end{itemize}
\end{lemma}

\noindent The cycles of gaps $\pgap(\pml{p})$ under this recursion form a discrete dynamic system.  

\vspace{0.08in}

A constellation of length $J$ is a sequence of $J$ consecutive gaps.
From here forward, we will focus on constellations $s$ of length $J$ and treat gaps as the special case $J=1$.

Under the recursion, constellations of certain lengths and spans are propagated in specific predictable ways \cite{FBHSFU, FBHPatterns}
from $\pgap(\pml{p_k})$ into $\pgap(\pml{p_{k+1}})$.
Two important observations facilitate the formulation of the population models for constellations of gaps.

\begin{corollary}\label{Obs1}
In step R3 of the recursion $\pgap(p_k^\#) \fto \pgap(p_{k+1}^\#)$ each of the $\phi(p_k^{\#})$ possible fusions 
in $\pgap(p_k^\#)$ occurs exactly once.
\end{corollary}

This is a result of the Chinese Remainder Theorem in this context.

\begin{corollary}\label{Obs2}
In step R3 of the recursion the minimum distance between fusions is $2 p_{k+1}$.
\end{corollary}

There may be longer constellations $\tilde{s}$ that form $s$ under the fusions of step R3 in the recursion.  For example, if we consider
$s=2,10,2,10,2$ of length $5$, the constellation $\tilde{s}=2462642$ produces $s$ when the pairs of adjacent gaps $46$ and $64$ 
are fused into $10$'s.
We call these constellations {\em driving terms} for $s$.

For a driving term $\tilde{s}$ of length $j \ge J$ for $s$, the $J+1$ fusions that would eliminate the image from being a driving term for 
$s$ are the {\em boundary fusions} for $\tilde{s}$.  The $j-J$ fusions for which the image of $\tilde{s}$ is still a driving term
for $s$ are the {\em interior fusions} for $\tilde{s}$.  For $s$ itself all $J+1$ of the fusions are boundary fusions.

\begin{figure}[hbt]
\centering
\includegraphics[width=3.1in]{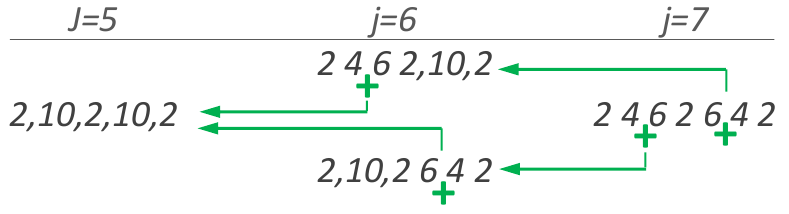}
\caption{\label{J5DynFig} For $s=2,10,2,10,2$ of length $J=5$, we show its admissible driving terms, and we mark its interior fusions
under the recursion.}
\end{figure}

If $|s| \le 2 p_{1}$, then the populations of $s$ and its driving terms in $\pgap(p_0^\#)$ form a Markov chain, and the associated
linear system has a beautifully simple eigenstructure  \cite{FBHSFU, FBHPatterns}, from which we can exactly model the populations in all
subsequent cycles of gaps $\pgap(p_k^\#)$.

\begin{figure}[hbt]
\centering
\includegraphics[width=4.25in]{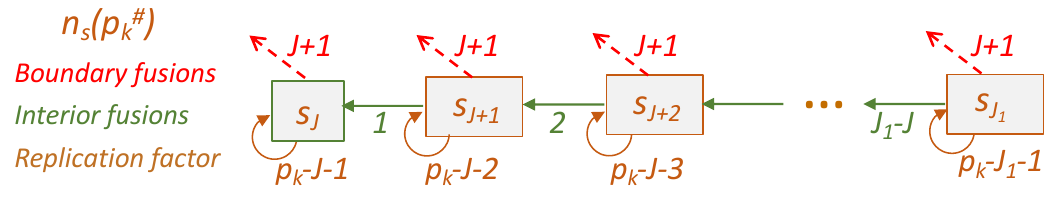}
\caption{\label{DynFig} Under the recursion, once $|s| \le 2 p_1$ then the populations of $s$ and its driving terms of length $j$ form
a Markov chain for all subsequent $p_k$.  In creating $\pgap(p_k^\#)$ the dynamics for all admissible constellations of length $J$ 
and span  $|s|< 2 p_k$ are identical.}
\end{figure}

\section{Admissible constellations}
A $k$-tuple is given in terms of generators, and a constellation is given in terms of the gaps between successive generators.  A constellation of length $J$ corresponds to a family of $k$-tuples with ${k = J+1}$.  

An {\em instance} of a constellation
$s$ is identified by specifying the initial generator $\gamma_0$.  Once $\gamma_0$ is set, the other generators in the $k$-tuple are 
determined by ${\gamma_j = \gamma_{j-1} + g_j}$.

For a constellation $s$ of length $J$ the parameter $\nu_p(s)$ is the number of residue classes $\bmod p$ covered by any instance
of $s$, with $1 \le \nu_p(s) \le \min\set{p, J+1}$.

We call $r=\gamma_0 \bmod p$ an {\em admissible instance} of $s$ in $\pgap(\pml{p})$ iff
$\gamma_j \neq 0 \bmod p$ for all $0 \le j \le J$.
The boundary fusions for $s$ are preserved for this $\gamma_0$ in $\pgap(\pml{p})$.  An instance of $s$
could be $s$ itself or a driving term.

There are $p-\nu_p(s)$ admissible instances $\gamma_0$ for $s$, modulo $p$. 
Define $\Upsilon$ to be the set of admissible values for the initial generator $ \gamma \bmod p$.
$$\Upsilon_s(p) = \set{ r = \gamma \bmod p \st \gamma_0=\gamma \; {\rm is~admissible~for} \; s}$$
The cardinality of $\Upsilon_s(p)$ is $p-\nu_p(s)$.

\begin{table}[h]
\centering
\begin{tabular}{rcc} \hline
\multicolumn{3}{c}{Admissibility for $s=2,10,2,10,2$.} \\
\multicolumn{1}{c}{$p$} & $\nu_p(s)$ & $\Upsilon_s(p) $ \\ \hline
$5$ & $4$ & $\set{2}$ \\
$7$ & $4$ & $\set{1,3,6}$ \\
$11$ & $5$ & $\set{1,2,3,4,5,6}$ \\
$13$ & $5$ & $\set{3,4,5,6,7,8,9,10}$ \\
$17$ & $6$ & $\set{1,2,4,6,7,9,11,12,13,14,16}$ \\
\end{tabular}
\end{table}
\noindent For $s=2,10,2,10,2$ and $p \ge 17$, $|s|<2p$ so all fusions occur in separate images under step R2, and ${\nu_p(s)=6}$.

We say that a constellation $s$ is {\em admissible for $p$} iff $\nu_p(s) < p$, or equivalently iff $\Upsilon_s(p)$ is non-empty.  
By default the constellation $s$ is admissible for all primes $p > J+1$.

\begin{lemma}\label{Lemdiv}
Let $s = g_1 g_2 \ldots g_J$ be a constellation of length $J$.
In step R3 of creating $\pgap(\pml{p})$, the fusions at $\gamma_i$ and $\gamma_j$ occur in the same image of $s$
iff $p$ divides the span $\left| \gamma_i - \gamma_j \right| = \sum_{i+1}^j g_{\iota}$.
\end{lemma}

{\bf Proof.}
The fusions in step R3 occur when we remove multiples of $p$ from the list of remaining candidate primes.
The fusion at $\gamma$ occurs when ${\gamma \bmod p = 0}$.  The residues $\bmod \, p$ for the $\gamma_i$ are set once
$\gamma_0$ is fixed.  If $\gamma_i \bmod p = 0$ in an image of $s$ in step R2, and $p$ divides ${\left| \gamma_i - \gamma_j \right|}$,
then $\gamma_j \bmod p = 0$ in this same image of $s$.
$\square$

\vspace{0.125in}

Finally, a constellation $s$ is {\em admissible} iff $s$ is admissible for all primes $p$.

\vspace{0.1in}

Let $Q(s)$ be the product of odd primes that divide a span between boundary fusions in $s$.  For the example $s=2,10,2,10,2$ this product
${Q(s)=3\cdot 5 \cdot 7 \cdot 11 \cdot 13}$.

\begin{lemma}\label{Lemnu}
Let $s$ be an admissible constellation of length $J$.  
\begin{enumerate}
\item[a.] If $p \le J+1$, then $p \mid Q(s)$. 
\item[b.] If $p > J+1$, then $\nu_p(s) < J+1$ iff $p \mid Q(s)$.
\end{enumerate}
\end{lemma}

{\bf Proof.} 
This is a remark on the top of p.62 in \cite{HL}.
For (a), when $p \le J+1$ there are at least as many boundary fusions as residue classes. Since $s$ is admissible, there must be multiple fusions
in some residue classes, and by Lemma~\ref{Lemdiv} this requires ${p \mid Q(s)}$.

When $p > J+1$ there are more residue classes than boundary fusions, and the admissibility of $s$ is automatic.  
In this case, we note that 
$$ p - \nu_p(s) > p-J-1 \gap {\rm iff} \gap p \mid Q(s).$$
If all $J+1$ boundary fusions occur in separate images of $s$, then $\nu_p(s)=J+1$.  Otherwise some boundary fusions occur in the
same image of $s$.  Lemma~\ref{Lemdiv} requires $p \mid Q(s)$, and $\nu_p(s) < J+1$.
$\square$

\section{Admissible instances across Eratosthenes sieve}
Under the recursion, Eratosthenes sieve systematically produces all admissible instances of every admissible constellation.
The admissible instances of $s$ include its driving terms.

We begin with a cycle of gaps $\pgap(\pml{p_0})$ for a $p_0$ that is manageable.  If needed, we could begin with the trivial starting points
$p_0=3$ or $5$, but we often start with ${7 \le p_0 \le 37}$.

\begin{figure}[hbt]
\centering
\includegraphics[width=4in]{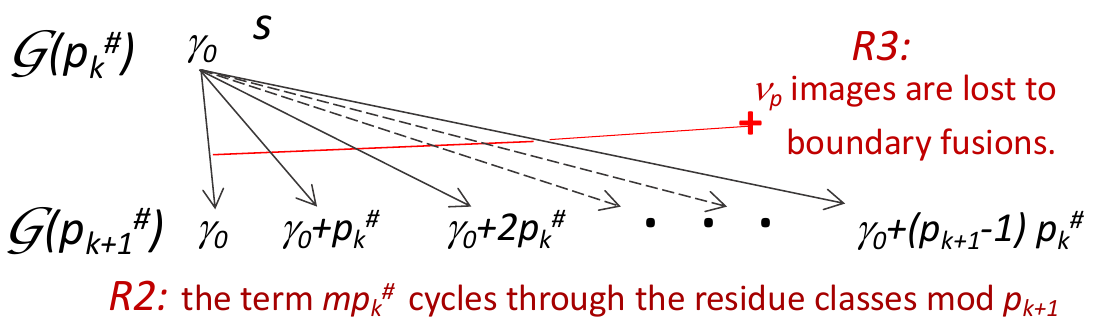}
\caption{\label{RFig} Under the recursion each instance $\gamma_0$ of $s$ in $\pgap(\pml{p_k})$ is replicated $p_{k+1}$ times
at offsets of $\pml{p_k}$.  Then $\nu_{p_{k+1}}$ of these images are lost to boundary fusions and $p_{k+1}-\nu_{p_{k+1}}$ survive
in $\pgap(\pml{p_{k+1}})$. }
\end{figure}

Eratosthenes sieve methodically eliminates all inadmissible instances of a constellation $s$, and the sieve methodically produces all admissible 
instances of $s$.  If a constellation $s$ is not admissible for the prime $p_{k+1}$, then every instance of $s$ and its driving terms is eliminated by
the recursion from ${\pgap(\pml{p_k}) \fto \pgap(\pml{p_{k+1}})}$.

\begin{lemma}\label{LemAdmIter}
Let $\gamma_0$ be the initial generator for an admissible instance of $s$ or a driving term $\tilde{s}$ in $\pgap(p_0^{\#})$.
Let
$$ \gamma_{0,k} = \gamma_0 \; + \; m_1 \cdot p_0^{\#} \; + \; m_2 \cdot p_1^{\#}  + \cdots + m_k \cdot p_{k-1}^{\#}$$
be an admissible instance of $s$ in $\pgap(\pml{p_k})$, with $0 \le m_j < p_j$.  
Then in step R2 of the recursion
${\pgap(\pml{p_k}) \fto \pgap(\pml{p_{k+1}})}$, this instance $\gamma_{0,k}$ of $s$ has $p_{k+1}$ images, with initial generators
$$ \left(\gamma_0 \; + m_1 \cdot \pml{p_0} \; +  \cdots + m_k \cdot p_{k-1}^{\#}\right) + m \cdot p_k^{\#} $$
for $m= 0,1,\ldots,p_{k+1}-1$.  These images cycle through all of the residue classes $\bmod \; p_{k+1}$.
As such, ${p_{k+1}-\nu_{p_{k+1}}(s)}$ of these survive in $\pgap(p_{k+1}^{\#})$, and $\nu_{p_{k+1}}(s)$ are eliminated by boundary fusions.
\end{lemma}

{\bf Proof.}
The replication $\gamma_{0,k}+m\cdot p_k^{\#}$ of this instance of $s$ is induced by step R2 of the recursion.
We have to show that these images $\gamma_{0,k}+m \cdot p_k^{\#}$ cycle through the residue classes $\bmod p_{k+1}$,
and that $p_{k+1}-\nu_{p_{k+1}}(s)$ of these survive the fusions.

Since $\pml{p_k} \neq 0 \bmod p_{k+1}$, the initial generators ${\gamma_{0,k} + m \cdot p_{k}^{\#}}$ cycle once through the residue classes
modulo $p_{k+1}$, for $0 \le m < p_{k+1}$.

For each $m$, $\gamma_{0,k}+m\cdot p_k^{\#}$ survives the fusions in step R3 iff
$$(\gamma_{0,k}+ m\cdot p_k^{\#} ) \bmod p_{k+1} \in \Upsilon_s(p_{k+1}).$$  
These fusions correspond to removing
multiples of $p_{k+1}$ from $\pgap(p_{k+1}^{\#})$, and each such multiple is ${0 \bmod p_{k+1}}$.
The cardinality of $\Upsilon_s(p_{k+1})$ is ${p_{k+1}-\nu_{p_{k+1}}(s)}$, and we have our result.
$\square$

\begin{figure}[hbt]
\centering
\includegraphics[width=3.1in]{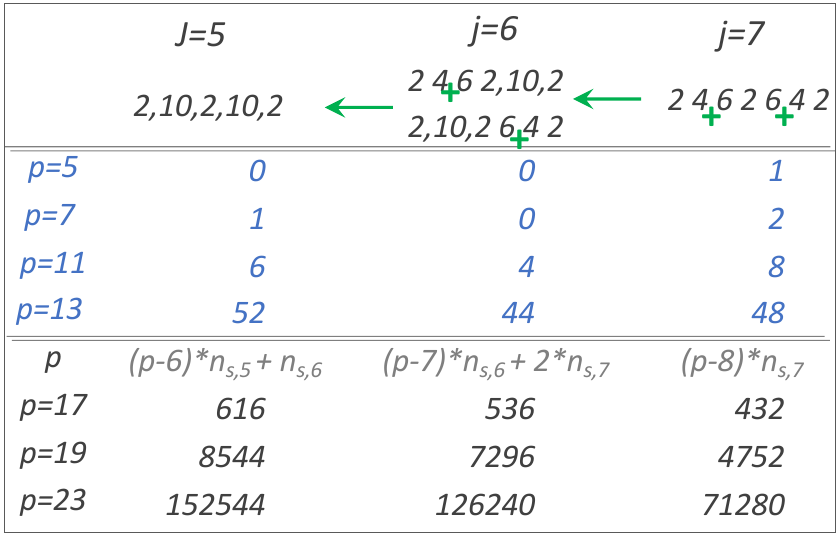}
\caption{\label{segFig} The populations $n_{s,j}$ of $s=2,10,2,10,2$ and its driving terms across the first few primes.
From $p=13$ on, the condition ${|s|<2p_1}$ is satisfied and the linear system holds. }
\end{figure}


An immediate result is Theorem~\ref{ThmAdmS}.

{\bf Proof of Theorem~\ref{ThmAdmS}.}
Let $s$ be a constellation with an instance at $\gamma_0$ in $\pgap(\pml{p_0})$.  The instance may be a driving term.
Suppose $s$ is admissible for all $p$ for $p_0 \le p \le p_k$.

For $i = 1,\cdots,k$, let $r_i \in \Upsilon_s(p_i)$ with ${0 \le r_i < p_i}$.  By Lemma~\ref{LemAdmIter} there exists a
unique $m_i$ with ${0 \le m_i < p_i}$ such that for 
$$\gamma_{0,i} = \gamma_{0,i-1}+m_i \cdot p_{i-1}^{\#}$$
$\gamma_{0,i} \bmod p_i = r_i$.  Since ${p_{i-1}^{\#} \bmod p = 0}$ for all ${p \le p_{i-1}}$, 
$$\gamma_{0,i} \bmod p_j = r_j \; \forall \; 0 \le j \le i.$$
For each ordered tuple of admissible $r_i$'s, there is a unique corresponding instance $\gamma_{0,k}$ in $\pgap(p_k^{\#})$.
$\square$

\vspace{0.1in}

By Theorem~\ref{ThmAdmS} all admissible instances of every admissible constellation arise in $\pgap(\pml{p})$.
Corollary~\ref{CorNsj} will show that the dynamics across $\pgap(\pml{p})$
treat all instances fairly by length $J$. 

The proofs of Lemma~\ref{LemAdmIter} and Theorem~\ref{ThmAdmS} yield the following corollary as well.

\begin{corollary}
Let $\gamma_{0,1}$ for $s_1$ and $\gamma_{0,2}$ for $s_2$ be instances of any two constellations in $\pgap(p_k^\#)$, of possibly different lengths.
Then the images of these instances of $s_1$ and $s_2$ under step R2 of the recursion cycle through the residue classes $\bmod \, p_{k+1}$
in the same order.  The starting points $\gamma_{0,1} \bmod p_{k+1}$ and $\gamma_{0,2} \bmod p_{k+1}$ may differ.
\end{corollary}

\section{Markov chain evolution.}

The first part of Corollary~\ref{CorNsj} follows directly from Theorem~\ref{ThmAdmS} and its proof.
Let $s$ be an admissible constellation of length $J$.  Then the count of admissible instances of $s$ in 
$\pgap(\pml{p})$ is
$$ \sum_{j \ge J} n_{s,j}(\pml{p}) = \prod_{q \le p} (q-\nu_q(s)).$$

This count includes the driving terms for $s$.
What is left to show is that under the dynamics of the Markov chain depicted in Figure~\ref{DynFig}, 
the populations shift from the longer driving terms to copies of $s$ itself, and to calculate the asymptotic
value for the relative population weight $w_{s,J}(\pml{p})$.

\begin{figure}[hbt]
\centering
\includegraphics[width=4.25in]{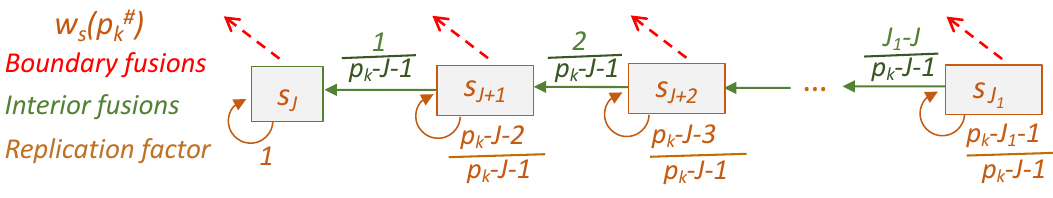}
\caption{\label{wJdynFig} When ${|s|<2p_1}$, the fusions all occur in separate images of $s$ under the recursion.
The system for the relative populations $w_{s,j}(\pml{p})$ then forms a Markov chain with the weights shown. }
\end{figure}

The population $n_{s,J}(\pml{p})$ is the population of the constellation $s$ in the cycle of gaps $\pgap(\pml{p})$.
These populations grow super-exponentially, dominated by factors of $(q-J-1)$.
We define the relative population $w_{s,J}(\pml{p})$  by
\begin{equation}
w_{s,J}(\pml{p}) = n_{s,J}(\pml{p}) \gap / \prod_{J+1 < q \le p} (q-J-1). \label{Eqwn}
\end{equation}
This Equation~(\ref{Eqwn}) isolates the dominant factors ${(q-J-1)}$.  If we take initial conditions in a cycle $\pgap(p_0^\#)$
for which ${|s| < 2p_1}$, then the relative population $w_{s,J}(\pml{p})$ 
is an iterative linear system as depicted in Figure~\ref{wJdynFig}.

The transfer matrix $M_J(\pml{p})$ at each stage is a banded matrix with diagonal elements $\frac{p-J-j}{p-J-1}$ and
super-diagonal elements $\frac{j}{p-J-1}$.  We then have
$$
w_s(p_{k}^\#)  = M_J(p_k) \cdot w_s(p_{k-1}^\#) \;  = \; M_J^k \cdot w_s(p_{0}^\#),
$$
using the notation
\begin{eqnarray*}
M_J^k & = & M_J(p_k)  \cdot M_J(p_{k-1}) \cdots  M_J(p_1) \\
    & = & R \cdot \Lambda_J^k \cdot L^T.
 \end{eqnarray*}
The matrix of right eigenvectors $R$ is an upper triangular Pascal matrix of alternating sign.  The matrix of left eigenvectors
$L^T$ is an upper triangular Pascal matrix.

The matrix of eigenvalues $\Lambda^k$ is
\begin{eqnarray*}
\Lambda_J^k & = & {\rm diag} \left(1, \gap \prod_{p_1}^{p_k}\frac{p-J-2}{p-J-1}, \gap \prod_{p_1}^{p_k}\frac{p-J-3}{p-J-1}, \gap \ldots, 
  \gap \prod_{p_1}^{p_k} \frac{p-J_1-1}{p-J-1} \right) \\
 & = & {\rm diag}\left( 1,  \gap \lambda=a_2^k, \gap a_3^k, \ldots , \gap a_{J_1-J+1}^k \right)
\end{eqnarray*}
where $J_1 \ge J$ is the length of the longest admissible driving term for $s$.

From the top row of
\begin{equation}
 w_s(p_k^\#) = M_J^k \cdot w_s(p_0^\#), \label{Eqw}
\end{equation}
we extract an equation for the relative population of $s$ itself
$$
w_{s,J}(p_k^\#) = \ell_1 \; - \; \ell_2 a_2^k \; + \; \ell_3 a_3^k - \cdots + (-1)^{i+1} \ell_{i} a_{i}^k + \cdots
$$
with coefficients $\ell_i = L_i^T\cdot w_s(p_0^\#)$.

\vspace{0.1in}

We note here that the system dynamics, e.g. the eigenvectors and eigenvalues, depend only on the primes and the length $J$.  
They do not vary by the specific
constellation of length $J$.  The various admissible constellations of length $J$ differ only in their initial conditions $w_s(p_0^\#)$, including the
length of their longest admissible driving term.

\vspace{0.1in}

{\bf Completing the proof of Corollary~\ref{CorNsj}.}
Let $s$ be an admissible constellation of length $J$, and let $Q(s)$ be the product of all odd primes that divide a span 
between boundary fusions in $s$.

The dominant eigenvalue in $\Lambda_J^k$ is $1$.  The other eigenvalues
$$
a_i^k \; = \; \prod_{p_1}^{p_k} \frac{p-J-i}{p-J-1}
$$
are all positive, and all converge to $0$ as $p_k \fto \infty$.  Setting $\lambda = a_2^k$ as our system parameter, we have the
approximation $a_i^k \approx \lambda^{i-1}$.  We can thus view the evolution of the model $w_{s,J}(p_k^\#)$ through the
evolution of the parameter $\lambda=a_2^k$.

Merten's Third Theorem tells us that 
$$\lambda = \mathcal{O}\left(\frac{1}{\ln p_k}\right) \fto 0$$
as $p_k \fto \infty$.

As $p_k \fto \infty$, we have the asymptotic value ${w_{s,J}(\infty) = \ell_1 = L_1^T\cdot w_s(p_0^\#)}$ for any prime
$p_0$ such that $|s| < 2p_1$.
The first left eigenvector is all $1$'s, so 
\begin{eqnarray*} 
w_{s,J}(\infty) & = & \ell_1 \cdot w_s(p_0^\#) = \sum_J^{J_1} w_{s,j}(p_0^\#) \\
  & = & \left( \sum_J^{J_1} n_{s,j}(p_0^\#) \right) / \prod_{J+1 < q \le p_0} (q-J-1) \\
  & = &  \prod_{q \le J+1 } (q-\nu_q(s)) \cdot  \prod_{J+1 < q \le p_0} \frac{q-\nu_q(s)}{q-J-1}  \\
  & = &  \prod_{q \le J+1 } (q-\nu_q(s)) \cdot  \prod_{\substack{J+1 < q \\  q | Q(s)}} \frac{q-\nu_q(s)}{q-J-1}  
\end{eqnarray*}
The last equality follows from Lemma~\ref{Lemnu}, and we have Corollary~\ref{CorNsj}.
$\square$

\begin{figure}[hbt]
\centering
\includegraphics[width=4.5in]{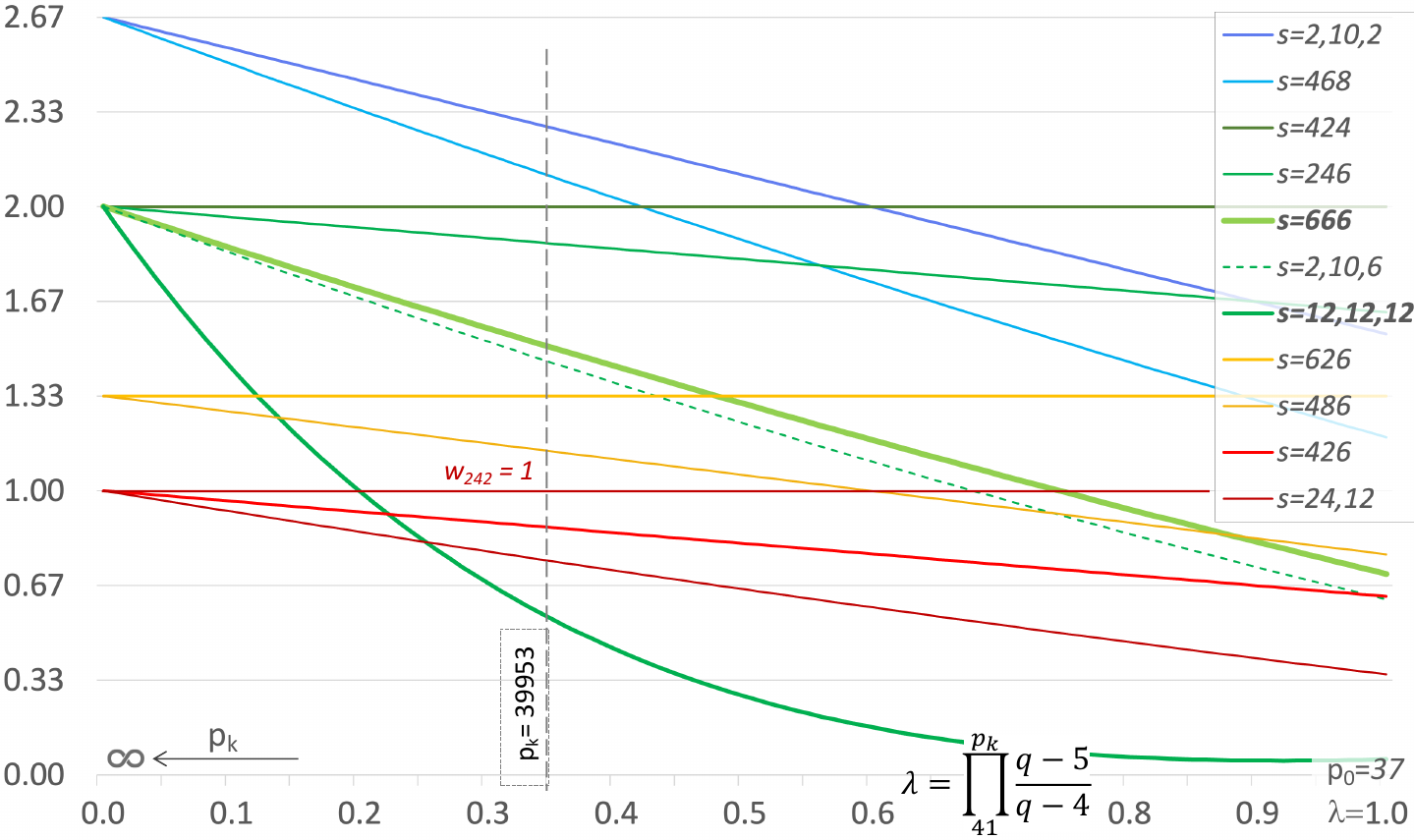}
\caption{\label{wJ3Fig} The relative populations $w_{s,3}$ shown for a few constellations of length $J=3$ with
$p_0=37$.  The relative populations of these constellations in $\pgap(37^\#)$ are on the right, at $\lambda=1$.
Their asymptotic values are shown on the left, at $\lambda=0$.  The section marked at $p=39953$ corresponds to
the sampled counts displayed in Figure~\ref{DelH3Fig}.}
\end{figure}

\vspace{0.125in}

A strong form of the $k$-tuple conjecture holds for Eratosthenes sieve.  Every admissible instance of every admissible constellation
arises in Eratosthenes sieve.  Moreover, the populations of all admissible constellations of length $J$ ultimately grow
at the same rate 
$$ n_{s,J}( p^\#) = \Theta\left( \prod_{q > J+1}^p (q-J-1) \right)$$
That is, for any two admissible constellations $s_1$ and $s_2$ of length $J$, their asymptotic relative populations are in a constant
proportion to each other.
$$ \lim_{p \rightarrow \infty} \frac{n_{s_1,J}(p^\#)}{n_{s_2,J}(p^\#)}
 \; = \; \lim_{p \rightarrow \infty} \frac{w_{s_1,J}(p^\#)}{w_{s_2,J}(p^\#)} \; = \; \frac{w_{s_1,J}(\infty)}{w_{s_2,J}(\infty)}
 $$
 From Corollary~\ref{CorNsj} the value of $w_{s,J}(\infty)$ can be calculated over the factors of $Q(s)$.
 
Further, every shorter admissible constellation eventually outnumbers every longer admissible constellation.
Let $s_1$ be an admissible constellation of length $J_1$, and $s_2$ be admissible of length $J_2$, with $J_1 < J_2$.
Then
$$ \lim_{p \rightarrow \infty} \frac{n_{s_2,J_2}(p^\#)}{n_{s_1,J_1}(p^\#)} 
 \; = \; c \cdot  \lim_{p \rightarrow \infty} \prod_{J_2 +1}^p \frac{q-J_2-1}{q-J_1-1}   \; = \; 0,
$$
where the constant $c$ is the ratio of the $w_{s,J}(\infty)$ and the factors $(q-J_1-1)$ up through $J_2$. 
This constant does not stop the slow decay of the product to $0$.

{\em Every} shorter admissible constellation (no matter how rare) eventually outnumbers {\em every} longer admissible constellation
(no matter how common) in this strong sense -- that the ratio of their populations in the cycles of gaps $\pgap(\pml{p})$ goes to $0$.
The evolutions of these populations take place over massive scales.  By Merten's Third Theorem the convergence goes
roughly as $1 / \ln(p)$.

One consequence is that it is meaningless to compare the relative populations $w_{s,J}(p^\#)$ for constellations of different lengths.
The population of the longer constellation in $\pgap(\pml{p})$ will always be asymptotically insignificant compared to the population of the
shorter constellation.

\section{Constellations of primes}
We have been studying the populations of constellations within the cycles of gaps $\pgap(\pml{p})$.  
How is this related to the populations of constellations among the prime numbers?

The constellations among primes are the constellations
toward the front of the cycles of gaps that survive further fusions.

\begin{figure*}[hbt]
\centering
\includegraphics[width=\textwidth]{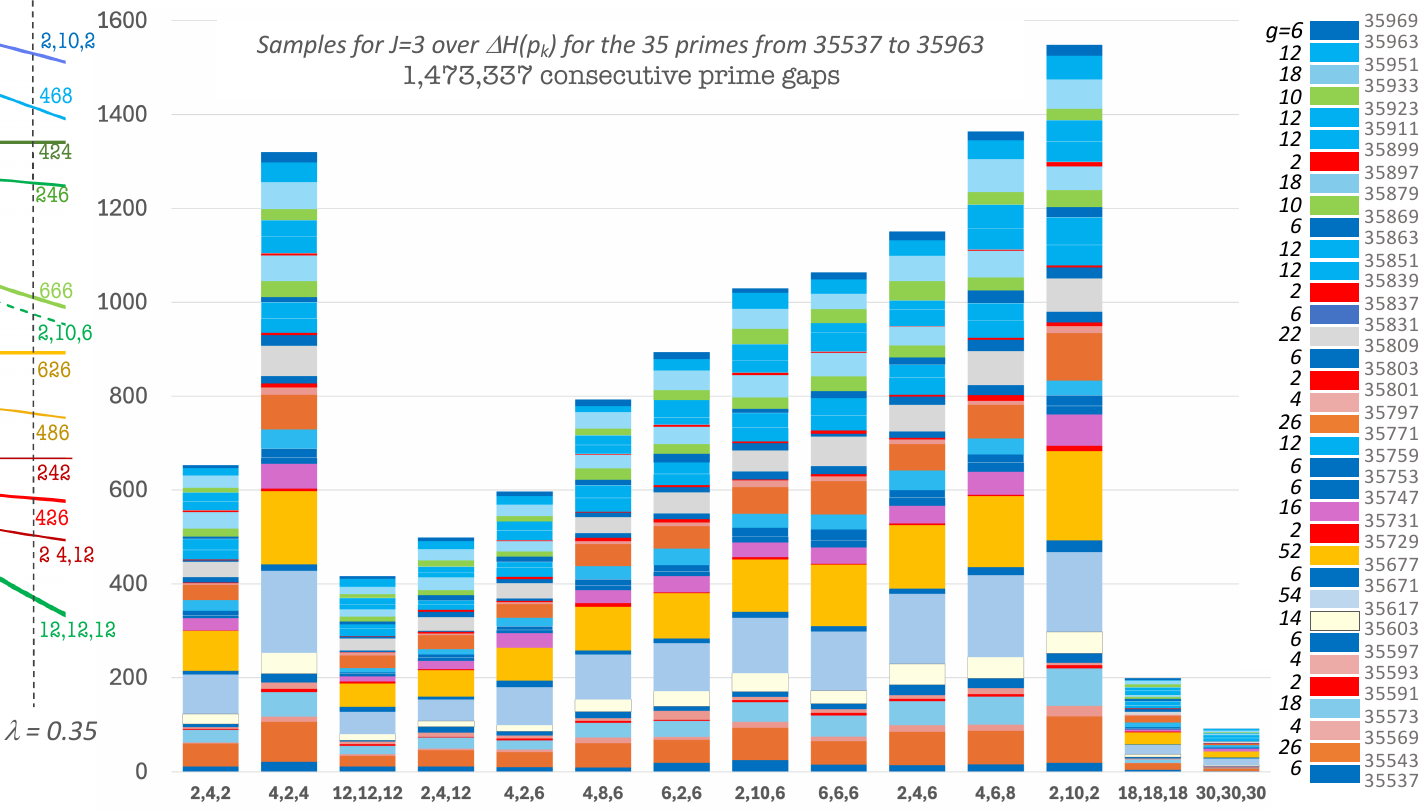}
\caption{\label{DelH3Fig} Counts of the occurrences of a few constellations of length $J=3$, among primes, in 35
intervals of survival $\Delta H(p)$ for the primes ${35537 \le p < 35969}$.  These intervals cover the $1\,473\,337$ primes between
$1\,262\,878\,369$ and $1\,293\,768\,961$.  The counts agree to first order with the section shown on the left, 
of $w_{s,3}(p^\#)$ from Figure~\ref{wJ3Fig} around $\lambda \approx 0.351$.  The legend at right shows the size of gap
$g=p_{k+1}-p_k$ for successive intervals of survival $\Delta H(p_k)$ and $\Delta H(p_{k+1})$.}
\end{figure*}

In the recursion from ${\pgap(p_{k}^\#) \rightarrow \pgap(p_{k+1}^\#)}$, the first fusion $g_1+g_2$ confirms $p_{k+1}$ as prime,
and the next fusion is at $p_{k+1}^2$.  We call $p_{k+1}^2$ the {\em horizon for survival} for the cycle $\pgap(p_k^\#)$.
All of the constellations between $p_{k+1}$ and $p_{k+1}^2$ will survive the sieve and be confirmed as constellations
among primes.

On the other hand, all of the constellations up through $p_{k}^2$ were confirmed in previous stages of the sieve.
So the interval that we would expect to best reflect the relative populations of constellations in $\pgap(p_k^\#)$
is the {\em interval of survival} from $p_{k}^2$ to $p_{k+1}^2$.  We denote this interval $\Delta H(p_k)$.
$$ \Delta H(p_k) \gap = \gap \left[ p_k^2, \; p_{k+1}^2 \right]$$

Our working hypothesis underlying our estimates for survival is the following:

\vspace{0.15in}

\noindent{\em Hypothesis on Survival:} Let $s$ be an admissible constellation of length $J$.
For $p$ much larger than $J$, the instances of $s$ are almost uniformly distributed in $\pgap(p^\#)$.

\vspace{0.15in}

The intuition behind this hypothesis is that under the recursion ${\pgap(p_k^\#)\rightarrow \pgap(p_{k+1}^\#)}$,
the $p_{k+1}$ images of each instance created in step R2 are in fact uniformly distributed at distances $p_k^\#$.
The fusions in step R3 do remove $J+1$ of these images, but for $p >> J$ these fusions will be a second order effect.
The fusions also occur symmetrically across $\pgap(\pml{p})$.  The cycle $\pgap(p^\#)$ is symmetric:
in $\pgap(p^\#)$, we have $ g_i \; = \; g_{\phi(p^\#)-i}$.

There is no part of the discrete dynamic system that is biased for or against any particular admissible constellation.
All admissible constellations of the same length grow asymptotically at the same rate, under the same dynamics.  
The relative populations already
account for the evolutions of different admissible constellations of the same length.

Consequently, as $p$ gets large we would expect the counts of constellations in the interval of survival $\Delta H(p)$
to reflect the relative populations $w_{s,J}(p^\#)$ at this stage of the sieve, to first order.
Figure~\ref{DelH3Fig} compares the actual populations of a few constellations of length $J=3$ across the intervals of
survival from $\Delta H(35537)$ to $\Delta H(35963)$.

\vspace{0.1in}

These models agree with Hardy and Littlewood  \cite{HL} in their conjecture ``Theorem X1".  In equation (5.664) they give
the weight
$$G(s) \; = \; \prod_{p \ge 2} \left( \left(\frac{p}{p-1}\right)^J \cdot \frac{p-\nu_p(s)}{p-1}\right)
$$
for the relative occurrence of the admissible constellation $s$ and its driving terms.  Their ``Theorem X1" depends upon a conjecture, and
they offer it as a consequence if the conjecture holds.  Further, their estimates are for differences as well as gaps.  In our setting,
this means that their estimates include counts of $s$ and all its driving terms.

They rewrite the weight $G(s)$ as a product ${G(s)=C_{J+1} H(s)}$, with $H(s)$ given by their Equation (5.667).
$$
H(s) \; = \;  \prod_{p \le J+1} \left( \left(\frac{p}{p-1}\right)^J \cdot \frac{p-\nu_p(s)}{p-1}\right) 
  \; \cdot \; \prod_{\substack{p \mid Q(s) \\ p > J+1}} \left( \frac{p-\nu_p(s)}{p-J-1}\right)
$$

Our Theorem~\ref{ThmAdmS} and Corollary~\ref{CorNsj} given above align strongly with this form of Hardy and Littlewood's k-tuple conjecture.
 Although we take a different approach to the subject, our results align with their estimates.

Initially the populations may be severely skewed toward long driving terms for $s$, but as the system evolves, eventually the occurrences of
$s$ itself will swamp the total population of all driving terms of lengths $j > J$.  We see the start of this evolutionary shift in the 
data for ${s=2,10,2,10,2}$ in Figure~\ref{segFig}.

\subsection{Consecutive primes in arithmetic progression}
Runs of consecutive primes in arithmetic progression (CPAP) are a special topic under the k-tuple conjecture \cite{BGreen}.
A run of $J+1$ consecutive primes in arithmetic progression corresponds to a repetition of length $J$ of a single gap $g$.

The following lemmas are easily proved and are offered here as remarks about repetitions of gaps.

\begin{lemma}
If $s$ is an admissible repetition of the gap $g$ of length $J$, and $p$ is the largest prime such that ${p \le J+1}$,
then $p^\# \mid g$.
\end{lemma}

\begin{lemma}
Let $s$ be an admissible repetition of the gap $g$ of length $J$, and let $Q$ be the product of odd primes dividing $g$.
Then the asymptotic relative population of $s$ in $\pgap(\pml{p})$ is
$$ w_{s,J} (\infty) \; = \;   \frac{\phi(Q)}{\prod_{q \mid Q, q > J+1}(q-J-1)}$$
\end{lemma}

This is a fairly high value for the asymptotic relative population in $\pgap(p^\#)$ but not always the maximum.
Further, the evolution of the relative populations of constellations is very slow.  Figure~\ref{wJ3Fig} shows the graphs
of $w_{s,3}(\lambda)$ for a few constellations of length $J=3$.  The graphs begin on the righthand side, with $p_0=37$
at $\lambda=1$, and they evolve to the left with parameter ${\lambda=\prod_{p_1}^{p_k} \frac{q-5}{q-4}}$.

The relative populations of the repetitions $s=666$ and $s=12,12,12$ are included in this graph.  Both have asymptotic value
$w_{s,3}(\infty)=2$.  At the section marked for $p_k=35963$, the parameter $\lambda \approx 0.35$, and 
\begin{eqnarray*}
w_{666,3}(0.35) & = & 1.503 \\
w_{12,12,12; 3}(0.35) & = & 0.547
\end{eqnarray*}

\begin{figure}[hbt]
\centering
\includegraphics[width=4.5in]{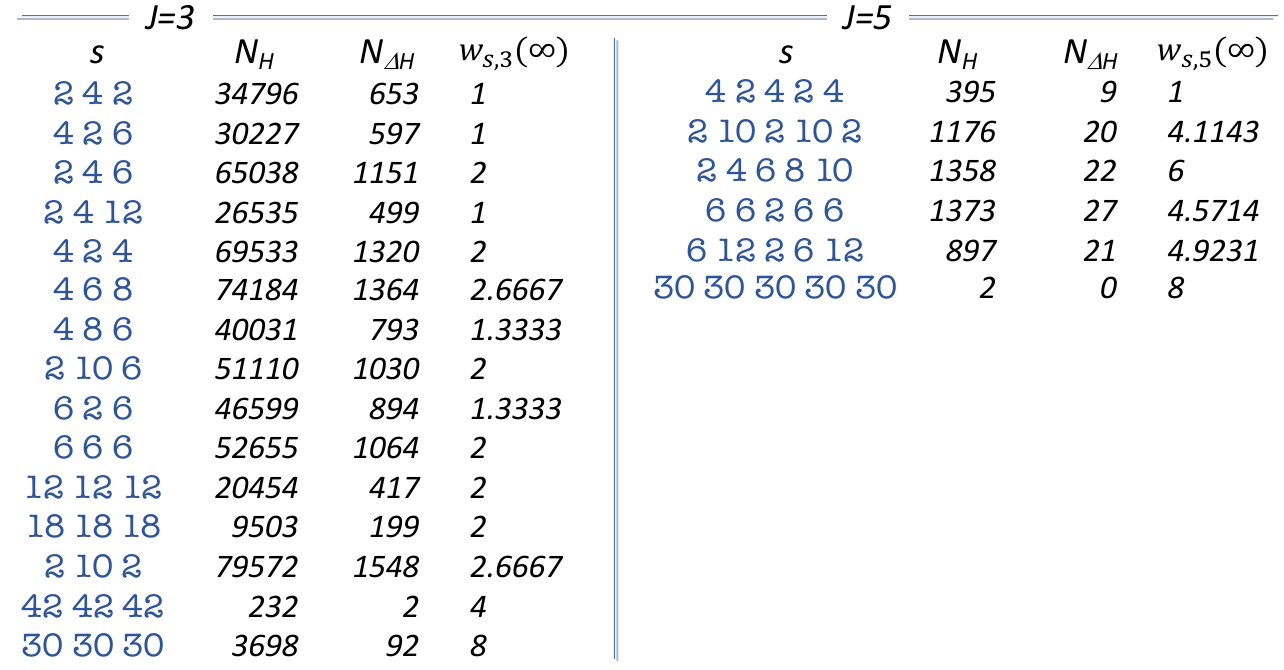}
\caption{\label{nHFig} We tabulate the populations of a selection of constellations for length ${J=3}$ and for $J=5$.
For each constellation we provide the total number of occurrences $N_H$ among primes from $29^2$ through $35969^2$,
and the number of occurrences $N_{\Delta H}$ across the $35$ intervals of survival 
from $\Delta H(35537)$ through $\Delta H(35963)$.  
The last column provides the asymptotic relative population, compared to other constellations of the same length.}
\end{figure}

Figure~\ref{nHFig} tabulates the total populations $N_H$ from $29^2$ to $35969^2$ for a few constellations of length $J=3$ and 
a few of length $J=5$.  The table also lists the populations $N_{\Delta H}$ for the constellations over the interval from 
$35537^2$ to $35969^2$, corresponding to the histogram in Figure~\ref{DelH3Fig}.  The asymptotic relative populations
$w_{s,J}(\infty)$ are shown in the last column.

We see that short repetitions (CPAP) are already cropping up.  The contrast between the counts and the asymptotic relative populations
highlight the evolution of the relative populations $w_{s,3}(p^\#)$, as shown in Figure~\ref{wJ3Fig}.
The section at $\lambda=0.35$ corresponds to the intervals represented in $N_{\Delta H}$.

The repetitions will become more frequent, compared to other constellations of the same length, but the sieve needs time to produce these.
CPAP are actually a relatively common occurrence but not until far far beyond the computational range.

In the meantime we can look for the driving terms for these repetitions.  If $s_{J+j}$ is a driving term for the admissible repetition
$s_J$, where $s_J$ has length $J$ and $s_{J+j}$ has length $J+j$, then surviving occurrences of the driving term $s_{J+j}$
will provide $J+j+1$ primes in arithmetic progression (AP).  

Green and Tao \cite{BGreen} cite known examples of primes in arithmetic progression of length $J+1=22$ and $J+1=23$.
We cast these examples in their primorial coordinates.
The shorter example with $J+1=22$ is a driving term for a repetition of the gap 
$$g  =  4609098694200  = 19^\# \cdot 20 \cdot 23 \cdot 1033$$ 
with 
\begin{eqnarray*}
\gamma_0 & = & 11410337850553  \\
 & = & 822463 + 3\cdot 19^\# + 19\cdot 23^\# +27\cdot 29^\#  +19\cdot 31^\# + 1\cdot 37^\#.
\end{eqnarray*}
The longer example with $J+1=23$ is a driving term for a repetition of the gap
$$g=44546738095860 = 23^\# \cdot 2 \cdot 99839$$ 
with 
\begin{eqnarray*}
\gamma_0 & = & 56211383760397 \\
 & = & 2164027 + 1\cdot 19^\# + 12\cdot 23^\# + 8\cdot 29^\#  + 21 \cdot 31^\# + 7\cdot 37^\#.
 \end{eqnarray*}

From our work above, we expect AP sequences to survive for $J=21$ and gap $g = 19^\# = 9699690$, and for $J=22$
and gap $g= 23^\# = 223092870$.  These will be driving terms for the repetitions, corresponding to CPAP.  These driving
terms exist in the cycles $\pgap(p^\#)$ and are relatively abundant among constellations of the same lengths.

\section{Conclusion}
We approach Eratosthenes sieve as a discrete dynamic system.  The objects of the system are the cycles of gaps $\pgap(p^\#)$
at each stage of the sieve.  There is a 3-step recursion that produces the next cycle of gaps from the current one.
$$ \pgap(p_k^\#) \; \fto \; \pgap(p_{k+1}^\#).$$
A lot of structure in $\pgap(p^\#)$ is preserved under this recursion.

Consequently if we take initial conditions from the cycle $\pgap(p_0)$, we can develop exact models for the relative populations
of all admissible constellations $s$ for which the span $|s| \le p_1$.  These models allow us to compare relative populations among
admissible constellations far beyond the computational horizon.  For computed samples of the populations of constellations, we can compare 
the samples to the values expected from the models.

We find that Eratosthenes sieve provides strong intuitive support for the $k$-tuple conjecture.  The sieve methodically 
produces every admissible instance of every admissible constellation.  Under the dynamics of the recursion, we find that
the populations of all admissible constellations of length $J$ ultimately grow at the same rate, $\prod (p-J-1)$.

This approach through Eratosthenes sieve aligns well with the estimates made by Hardy and Littlewood.

Code and data supporting the analysis above can be found at {\small https://github.com/fbholt}.

\bibliographystyle{integers}
\bibliography{primes2024}

\end{document}